\documentclass[12pt]{amsart}
\usepackage{amsthm}
\usepackage{amsmath}
\usepackage{amsfonts}
\usepackage{amssymb}
\input amssym.def
\usepackage{pstricks}
\newtheorem{theorem}{Theorem}

\newtheorem{corollary}[theorem]{Corollary}
\theoremstyle{definition}
\newtheorem{remark}[theorem]{Remark}
\newcommand{\M}{\mathcal{M}}
\newcommand{\Z}{\mathbb{Z}}
\newcommand{\R}{\mathbb{R}}
\newcommand{\T}{\mathcal{T}}
\author{B\l{}a\.zej Szepietowski}
\title{On the commutator length of a Dehn twist}
\address[]{Institute of Mathematics, Gda\'nsk University, Wita Stwosza 57,
80-952 Gda\'nsk, Poland} 
\email{blaszep@mat.ug.edu.pl}
\thanks{Supported by MNiSW grant N N201 366436}
\begin{document}

\begin{abstract}
We show that on a nonorientable surface of genus at least 7 any power of a Dehn twist is equal to a single commutator in the mapping class group and  the same is true, under additional assumptions, for the twist subgroup, and also for the extended mapping class group of an orientable surface of genus at least 3. 
\end{abstract}
\maketitle
\section{Introduction}
Let $S$ be a closed surface of genus $g$. If $S$ is nonorientable, then $g$ is the number of projective planes in a connected sum decomposition. The {\it mapping class group} $\M(S)$ of  $S$ is the group of isotopy classes of all, orientation preserving if $S$ is orientable, self-homeomorphisms of $S$.  For   orientable $S$, the {\it extended mapping class group} $\M^\diamond(S)$ is the group of isotopy classes of all self-homeomorphisms of $S$, including those reversing orientation. 
For notational convenience we define $\M^\diamond(S)$ to equal $\M(S)$ for nonorientable $S$. 

For a two-sided simple closed curve $c$ on  $S$ we denote by $t_c$ a Dehn twist about $c$. For a Dehn twist $t_c$ we always assume that $c$ does not bound a disc nor a M\"obius band, so that $t_c$ is a nontrivial element of $\M(S)$. It is well known that $\M(S)$ is generated by Dehn twists if $S$ is orientable. If $S$ is nonorientable, then the {\it twist subgroup} $\T(S)$ generated by all Dehn twists has index 2 in $\M(S)$ (cf. \cite{Lick2}).

For a group $G$  let $[G,G]$ denote the commutator subgroup generated by all commutators $[a,b]=aba^{-1}b^{-1}$. For $x\in[G,G]$ the {\it commutator length}
$cl_G(x)$ is the minimum number of factors needed to express $x$ as a product of commutators. The {\it stable commutator length} is the limit
\[scl_G(x)=\lim_{n\to\infty}\frac{cl_G(x^n)}{n}.\]
Recall that the first homology group $H_1(G;\Z)$ of $G$ is isomorphic to the quotient $G/[G,G]$.

For orientable $S$ and  $g\ge 3$ it is well known that $\M(S)$ is perfect, i.e. $\M(S)=[\M(S),\M(S)]$ (cf. \cite{Powell}), and for any Dehn twist $t_c$ we have  $cl_{\M(S)}(t_c)=2$ and $scl_{\M(S)}(t_c)>0$  \cite{BK,EndoKot,KorkOzb,KorkScl}. 
For nonorientable $S$ the groups $H_1(\M(S);\Z)$ and $H_1(\T(S);\Z)$ were computed by Korkmaz \cite{KorkHom} and Stukow \cite{Stukow}. In particular, if  $g\ge 7$ then we have
\[[\M(S),\M(S)]=[\T(S),\T(S)]=\T(S).\]
In this paper we prove the following. 

\begin{theorem}\label{extended}
Let $S$ be a closed orientable surface of genus $g\ge 3$ or
a closed nonorientable surface of genus $g\ge 7$.
 Then for every two-sided simple closed curve $c$ on $S$ and every $n\in\mathbb{Z}$,  
$t^n_c$ is equal to a single commutator of elements of $\M^\diamond(S)$. 
\end{theorem}

For even $n$ Theorem \ref{extended}, and hence also Corollary \ref{scl} below, follow immediately form the fact that $t_c$ is conjugate to its inverse in $\M^\diamond(S)$ (see Remark \ref{rem}), and were known already, at least for orientable $S$. See \cite[Remark 12]{Kot}.  

\begin{theorem}\label{twistsbgp}
Let $c$ be a two-sided simple closed curve on a closed nonorientable surface $S$
satisfying one of the following assumptions.  
\begin{itemize}
\item $c$ is separating and $g\ge 7$, or
\item $S\backslash c$ is connected and nonorientable and $g\ge 8$, or
\item $S\backslash c$ is connected and orientable, $g\ge 6$ and
$g\equiv 2\bmod 4$.
\end{itemize} 
Then for any $n\in\mathbb{Z}$,  
$t^n_c$ is equal to a single commutator of elements of $\T(S)$.
\end{theorem}

The following corollary is an immediate consequence of Theorems \ref{extended} and \ref{twistsbgp} and the definition of the stable commutator length.

\begin{corollary}\label{scl}
For $S$ and $c$ as in Theorem \ref{extended} or Theorem \ref{twistsbgp} we have
respectively $scl_{\M^\diamond(S)}(t_c)=0$ or $scl_{\T(S)}(t_c)=0$.
\end{corollary}

Our proof of Theorem \ref{twistsbgp} fails when $c$ is nonseparating and $g=7$, or $g=4k$ for $k\ge 2$ and $S\backslash c$ is orientable. We conjecture that also in these cases we have $cl_{\T(S)}(t^n_c)=1$ for any $n\in\mathbb{Z}$.

\section{Proofs}
\begin{figure}
\begin{center}
\input{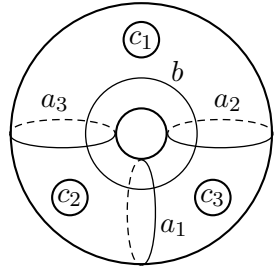} 
\end{center}
\caption{\label{torus}The torus $T$. 
}
\end{figure}

Consider a torus with three holes $T$. Let $c_1, c_2, c_3$ be its boundary curves and let $b, a_1, a_2, a_3$ be nonseparating simple closed curves 
in the interior of $T$, such that $a_1, a_2, a_3$ are pairwise disjoint, and $b$ intersects $a_i$ transversally at one point for $i=1,2,3$ (Figure \ref{torus}).
The right Dehn twists about these curves satisfy the following relations in the mapping class group of $T$.
\begin{itemize}
\item Twists about disjoint curves commute,
\item $t_bt_{a_i}t_b=t_{a_i}t_bt_{a_i}$ for $i=1,2,3$,
\item $(t_bt_{a_1}t_{a_2}t_{a_3})^3=t_{c_1}t_{c_2}t_{c_3}$. 
\end{itemize}
The first two are the well known braid relations, the third is the star relation
discovered by Garvais \cite{Gerv}.
By using the braid relations we can rewrite the star relation in the following way.
\[
\begin{aligned}
t_{c_1}t_{c_2}t_{c_3}&=(t_bt_{a_1}t_{a_2}t_{a_3})(t_bt_{a_1}t_{a_2}t_{a_3})(t_bt_{a_1}t_{a_2}t_{a_3})\\
&=t_bt_{a_2}t_{a_3}(t_{a_1}t_bt_{a_1})t_{a_2}(t_{a_3}t_bt_{a_3})t_{a_1}t_{a_2}\\
&=t_bt_{a_2}t_{a_3}t_bt_{a_1}(t_bt_{a_2}t_b)t_{a_3}t_bt_{a_1}t_{a_2}\\
&=(t_bt_{a_2}t_{a_3}t_bt_{a_1}t_{a_2})(t_bt_{a_2}t_{a_3}t_bt_{a_1}t_{a_2})
\end{aligned}
\]
Since $t_{c_i}$ commute with all twists, for every $n\in\mathbb{Z}$ we have
\begin{equation}\label{rel1}
t^n_{c_1}=(t_bt_{a_2}t_{a_3}t_bt_{a_1}t_{a_2}t^{-1}_{c_2})^n(t^{-1}_{c_3}t_bt_{a_2}t_{a_3}t_bt_{a_1}t_{a_2})^n.
\end{equation}
There is a reflectional symmetry $r\colon T\to T$ such that
$r(b)=b$, $r(a_1)=a_1$, $r(a_2)=a_3$, $r(c_1)=c_1$, $r(c_2)=c_3$. Since $r$ is orientation reversing it conjugates right twists to left twists, and so we have
\[r(t_bt_{a_2}t_{a_3}t_bt_{a_1}t_{a_2}t^{-1}_{c_2})^{-n}r=(t^{-1}_{c_3}t_{a_3}t_{a_1}t_bt_{a_2}t_{a_3}t_{b})^n.\]
By the braid relations we have
\[
\begin{aligned}
&t^{-1}_{c_3}t_{a_3}t_{a_1}t_bt_{a_2}t_{a_3}t_{b}=t^{-1}_{c_3}t_{a_1}(t_{a_3}t_bt_{a_3})t_{a_2}t_b=
t^{-1}_{c_3}t_{a_1}t_bt_{a_3}(t_bt_{a_2}t_b)=\\
&=t^{-1}_{c_3}t_{a_1}t_bt_{a_3}t_{a_2}t_bt_{a_2}=t_{a_1}(t^{-1}_{c_3}t_bt_{a_2}t_{a_3}t_bt_{a_1}t_{a_2})t^{-1}_{a_1}.
\end{aligned}
\]
Thus
\[(t^{-1}_{c_3}t_bt_{a_2}t_{a_3}t_bt_{a_1}t_{a_2})^n=t^{-1}_{a_1}r(t_bt_{a_2}t_{a_3}t_bt_{a_1}t_{a_2}t^{-1}_{c_2})^{-n}rt_{a_1},\]
and by (\ref{rel1})
\begin{equation}\label{rel2}
t^n_{c_1}=[(t_bt_{a_2}t_{a_3}t_bt_{a_1}t_{a_2}t^{-1}_{c_2})^n,t^{-1}_{a_1}r].
\end{equation}

\medskip

\begin{figure}
\begin{center}
\input{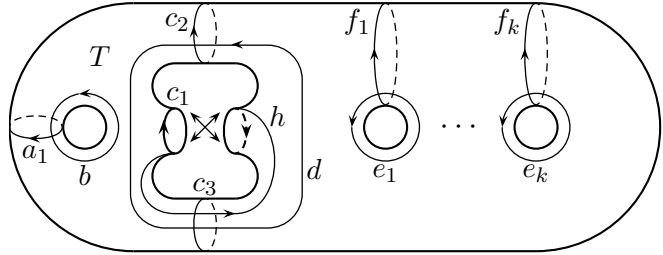} 
\end{center}
\caption{\label{f2}An embedding of the torus $T$ in a closed nonorientable surface $S$ such that $S\backslash c_1$ is orientable.
}
\end{figure}

{\it Proof of Theorem \ref{extended}.}
Let $S$ be a closed orientable surface of genus $g\ge 3$ or
a closed nonorientable surface of genus $g\ge 7$.
Consider the torus $T$ as embedded in  $S$ in such a way that the reflectional 
symmetry $r$ extends to $r\colon S\to S$. Then (\ref{rel2}) holds in $\M^\diamond(S)$. The embedding of $T$ in $S$ can be arranged in such a way that $c_1$ is 
nonseparating in $S$ or separating, bounding a subsurface of arbitrary topological type.
Moreover, if $S$ is nonorientable of even genus and $c_1$ is nonseparating  then $S\backslash c_1$ may be orientable or not (the former case is shown on Figure \ref{f2}). 
It follows that for every simple closed curve $c$ on $S$ there is a 
homeomorphism $h\colon S\to S$ such that $h(c)=c_1$, for appropriate embedding of $T$ in $S$. Thus $t^n_c=h^{-1}t^n_{c_1}h$ for $n\in\mathbb{Z}$ and, since a conjugate of a commutator is also a commutator, 
we have proved Theorem \ref{extended}.

\medskip

Recall that for a closed nonorientable surface $S$ of genus $g$, the homology group  $H_1(S;\R)$ is a real vector space of dimension $g-1$. For $f\in\M(S)$ let $f_\ast\colon H_1(S;\R)\to H_1(S;\R)$ be the induced automorphism. It turns out that the determinant homomorphism $f\mapsto\det f_\ast$ takes values in the group $\{-1,1\}$ and its kernel is the twist subgroup $\T(S)$ (cf. \cite{Lick1} and \cite[Corollary 6.3]{Stukow}).

\medskip

{\it Proof of Theorem \ref{twistsbgp}.} The idea of the proof is the same as for Theorem \ref{extended}. The only problem is that the involution $r$ may not be an element of $\T(S)$, in which case it has to be replaced by a different mapping class.

Suppose that $S$ and $c$ satisfy one of the assumptions of the theorem. 
Consider $T$ as embedded in $S$ in such a way that $c_1=c$, as in the proof of Theorem \ref{extended}. If $c$ is separating, then we may arrange that the component of $S\backslash c$ which does not contain $T$ is nonorientable of genus at least 2. If $c$ is separating  or nonseparting with $S\backslash c$ nonorientable, then $N=S\backslash T$ is a nonorientable surface of genus at least 2 and hence it supports a homeomorphism $h$ which is not a product of Dehn twists (we may take $h$ to be a crosscap slide or Y-homeomorphism introduced by Lickorish \cite{Lick1}). Since $h$ is equal to the identity on $T$, thus it commutes with twists about  $b$, $a_i$, $c_i$ for $i=1,2,3$. Now 
if the involution $r$, which is an extension of the reflectional symmetry of $T$, is in $\T(S)$, then 
$cl_{\T(S)}(t^n_c)=1$ by (\ref{rel2}). If $r\notin\T(S)$ then $rh\in\T(S)$, and 
$cl_{\T(S)}(t^n_c)=1$ by (\ref{rel2}) with $rh$ in the place of $r$.

Now suppose that $S\backslash c_1$ is orientable. Figure \ref{f2} shows $S$ as being obtained from an orientable surface $S'$ by identifying two boundary components. 
The homology classes of the curves $a_1$, $b$, $c_2$, $d$, $h$, $e_i$, $f_i$ for $1\le i\le k$, where $g=2(k+3)$, form a basis of $H_1(S;\R)$ (note that $[c_1]=0$ in $H_1(S;\R)$). 
Now we may take $r$ as being induced by a reflection of $S'$, so that
$r_\ast[a_1]=[a_1]$, $r_\ast[b]=-[b]$, $r_\ast[c_2]=[c_2]$, $r_\ast[d]=-[d]$,
$r_\ast[h]=[h]-[d]$, $r_\ast[e_i]=-[e_i]$ $r_\ast[f_i]=[f_i]$ for $1\le i\le k$.
We see that $\det r_\ast=(-1)^k$, which means that $r\in\T(S)$ if and only if $g\equiv 2\bmod 4$.

\medskip

\begin{remark}\label{rem}
Let $c$ be any two-sided simple closed curve on a surface $S$. There is a homeomorphism $h\colon S\to S$ preserving $c$ and reversing orientation of its neighbourhood. We have $t_c=ht^{-1}_ch^{-1}$ and 
\[t^{2n}_c=t^n_ct^n_c=t^n_cht^{-n}_ch^{-1}=[t^n_c,h]\]
for any $n\in\mathbb{Z}$. Thus any even power of any Dehn twist on any surface $S$ is equal to a single commutator of elements of $\M^\diamond(S)$. Moreover, if $S\backslash c$ is nonorientable of genus at least 2, then we may take $h\in\T(S)$ by composing it if necessary with a homeomorphism fixing $c$ which is not a product of Dehn twists. In particular, if $c$ is a nonseparating two-sided curve on a closed nonorientable surface of genus $7$ then $scl_{\T(S)}(t_c)=0$, which slightly improves Corollary \ref{scl}.
\end{remark}

\end{document}